\documentclass[pdflatex,sn-mathphys-num]{sn-jnl}

\usepackage{graphicx}
\usepackage{amsmath,amssymb}
\usepackage{booktabs}
\usepackage{tabularx,array}
\usepackage{xurl}

\let\orcidlogo\relax
\usepackage{orcidlink}

\newcolumntype{Y}{>{\raggedright\arraybackslash}X}

\AtBeginDocument{%
  \hypersetup{%
    pdftitle={The Stairs of Reconciliation: A Mathematical Tourist in Graz},%
    pdfauthor={Santiago Schnell},%
    pdfsubject={Mathematical Tourist column submission},%
    pdfkeywords={mathematical tourism, spiral staircase, intersecting circles, helicoid, minimal surfaces, mathematical biology}%
  }%
}

\begin{document}

\title[The Stairs of Reconciliation]{The Stairs of Reconciliation: A Mathematical Tourist in Graz}

\author*[1,2,3]{%
  \fnm{Santiago} \sur{Schnell}\,
  \orcidlink{0000-0002-9477-3914}%
}
\email{santiago.schnell@dartmouth.edu}

\affil*[1]{\orgdiv{Department of Mathematics}, \orgname{Dartmouth College},
\orgaddress{\city{Hanover}, \state{New Hampshire}, \postcode{03755}, \country{USA}}}

\affil[2]{\orgname{Geisel School of Medicine at Dartmouth},
\orgaddress{\city{Hanover}, \state{New Hampshire}, \postcode{03755}, \country{USA}}}

\abstract{Inside the Grazer Burg, two late-Gothic stone flights rise about distinct
spindles, overlap, share several treads, and separate again. Their plan is
governed not by a coaxial double helix but, to first approximation, by two
intersecting circles. This elementary geometry yields a model of recurrent
meeting and makes explicit the compatibility conditions that meeting requires.
It also leads to a second object that geometers call a double spiral staircase---
the helicoid---and to a useful distinction between resemblance and identity.
The staircase becomes a meditation on how paths, models, and disciplines can
meet without becoming the same.}

\keywords{mathematical tourism, spiral staircase, intersecting circles, helicoid, minimal surfaces, mathematical biology}

\maketitle

\section{An unplanned detour}
I had not come to Graz for a staircase. I had come for the 14th European
Conference on Mathematical and Theoretical Biology, the joint meeting of the
European Society for Mathematical and Theoretical Biology and the Society for
Mathematical Biology, held at the University of Graz in July 2026
\citep{ECMTB2026}. My wife had spent part of the conference exploring the city.
Knowing my susceptibility to geometric marvels, she led me, shortly before our
departure, along Hofgasse, through the courtyard of the Grazer Burg, and into a
narrow doorway.

\begin{figure}[t]
  \centering
  \includegraphics[width=0.96\linewidth]{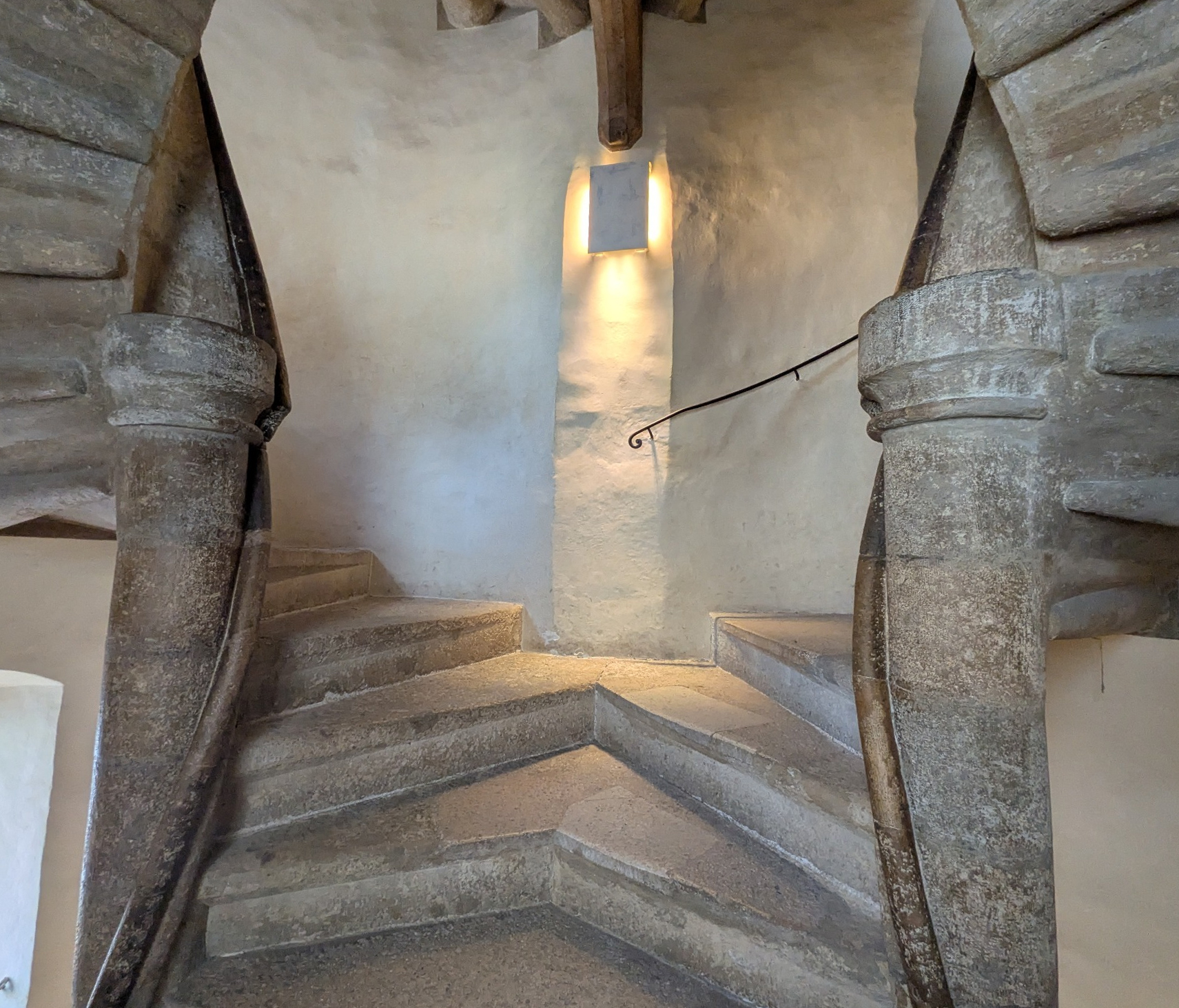}
  \caption{The Grazer \emph{Doppelwendeltreppe} viewed from a reunion zone. Two flights
  rise about distinct spindles, share several treads in the middle, and divide
  again. Photograph \copyright\ 2026 Mariana Schnell; used with permission}
  \label{fig:reunion}
\end{figure}

What waits inside is not a single spiral divided by a balustrade. It is a twin
stair. Two flights rise about two different spindles, one turning clockwise and
the other counter-clockwise. Between the spindles they approach, share several
stone treads, and divide again. Higher up, the pattern recurs. The effect is
Escherian not because the construction is impossible, but because it is exact:
the eye keeps losing track of which stair is one and which is the other
(Figure~\ref{fig:reunion}).

The stair was completed around 1499--1500 during the building campaign of
Emperor Maximilian~I (Holy Roman Emperor after 1508). Architectural descriptions call it a
\emph{Zwillingswendeltreppe}, a twin spiral staircase, and compare its horizontal
plan with a figure eight~\citep{Ginelli1986}. Graz residents also call it the
\emph{Vers\"ohnungsstiege}, the stairs of reconciliation (locally also the
\emph{Busserlstiege}, the kissing staircase~\citep{GrazStadtportal}): people who
separate on the two flights meet again where the stone becomes common
\citep{GrazTourismus,Castello2023}. Reconciliation is not the
erasure of difference; it is the restoration of relation. The Latin
\emph{reconciliare} means to restore or bring together again. Its relation to
\emph{concilium}, an assembly or council, lends the metaphor of common counsel an
apt resonance.

\section{A model with two centres}
The first mathematical task is to identify the object being modelled. A tempting
idealization places two oppositely handed helices on one cylinder about one axis.
That construction is elegant, but it is not the Graz stair. The two flights have
distinct centres. A model that suppresses tread width and replaces each flight by
a smooth centreline should therefore begin with two circles.

Let the spindle centres in plan be
\[
 C_{\rm L}=(-a,0),\qquad C_{\rm R}=(a,0),
\]
let both idealized walking lines have radius $R$, and assume
\[
 0<a<R.
\]
The two plan curves are
\begin{equation}
 (x+a)^2+y^2=R^2,
 \qquad
 (x-a)^2+y^2=R^2.
 \label{eq:plan-circles}
\end{equation}
Subtracting the equations gives $x=0$; substitution then gives the two
intersection points
\begin{equation}
 P_{\pm}=\left(0,\,\pm\sqrt{R^2-a^2}\right).
 \label{eq:plan-intersections}
\end{equation}
Thus two meetings occur precisely when the circular plans overlap. The limiting
cases are informative: if $a=R$ the circles are tangent and there is one meeting;
if $a>R$ they are disjoint and there is none. Opposite circulation by itself does
not force reconciliation.

\begin{figure}[t]
  \centering
  \includegraphics[width=0.92\linewidth]{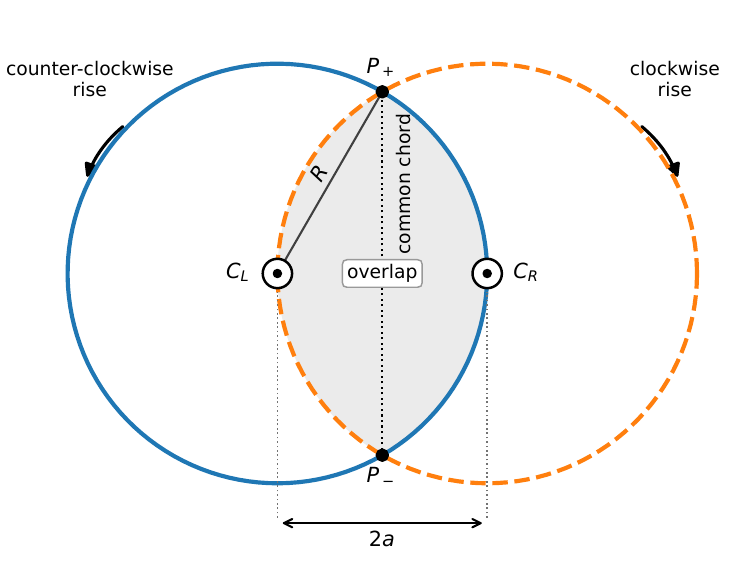}
  \caption{A symmetric two-centre plan model. The idealized flight centrelines
  are equal circles of radius $R$ with centres $C_L$ and $C_R$ separated by
  $2a$. For $0<a<R$ they meet at $P_+$ and $P_-$. The shaded lens marks the
  overlap in which the real stair can form shared treads. This is a qualitative
  model, not a measured survey}
  \label{fig:plan}
\end{figure}

To add height, assume mirror symmetry, a common vertical rise rate $p>0$, and a
common vertical phase. One convenient parametrization is
\begin{align}
 \Gamma_L(t)&=\bigl(-a+R\cos t,\;R\sin t,\;pt\bigr),
 \label{eq:left-flight}\\
 \Gamma_R(t)&=\bigl(a-R\cos t,\;R\sin t,\;pt\bigr).
 \label{eq:right-flight}
\end{align}
The left flight turns counter-clockwise around $C_L$. Relative to $C_R$, the
angular coordinate of the second is $\pi-t$, so the right flight turns clockwise
as $t$ increases. They have opposite handedness, but their reunion still depends
on the additional assumptions of overlap, equal rise rate, and synchronized
phase.

Put
\[
 \alpha=\arccos(a/R),\qquad
 b=\sqrt{R^2-a^2}=R\sin\alpha.
\]
Equations~\eqref{eq:left-flight} and \eqref{eq:right-flight} agree when
$\cos t=a/R$, hence
\begin{equation}
 \Gamma_L(2\pi k\pm\alpha)
 =\Gamma_R(2\pi k\pm\alpha)
 =\bigl(0,\,\pm b,\,p(2\pi k\pm\alpha)\bigr),
 \qquad k\in\mathbb Z.
 \label{eq:three-dimensional-meetings}
\end{equation}
The vertical separations between consecutive meetings alternate between
$2p\alpha$ and $2p(\pi-\alpha)$. They are equal only in the coaxial limiting case
$a=0$, for which $\alpha=\pi/2$. The real staircase has finite-width, discrete
stone treads, so the point intersections of the centreline model are thickened
into short common zones. Without a measured plan and section, the equations
should be read as a symmetric kinematic model, not as a dimensional
reconstruction.

A mathematical model earns its force by
stating what has been idealized, which relations are structural, and which
conclusions depend on chosen symmetries. Here the structural feature is not
simply ``two opposite spirals.'' It is two overlapping circular systems whose
vertical motions have been made compatible.

\section{Geometry in the workshop}
The two-centre model is not merely imposed after the fact. Anatol 
Ginelli~\citep{Ginelli1986} compared the plan with a figure eight. A later architectural
account says that the masons worked with compasses, brought two circles into
overlap, and treated their common chord as an axis of symmetry
\citep{WerkgruppeGraz}. The stone elements were cut with a constant rise ratio;
in the overlap, the symmetry cut determined the joint angles and made the two
flights assemblable as one construction. Much was resolved in a single cut:
each tread was shaped in the workshop as one stone joining its bearing, its
projecting winder, and its portion of the spindle, so that walking surface,
support, and core rose together as one repeated element.
A recent architectural-conservation study~\citep{Castello2023} likewise treats
the stair as a twin, rather than merely doubled, system and analyses its
historical, typological, and constructive logic.

\begin{figure}[t]
  \centering
  \includegraphics[width=0.75\linewidth]{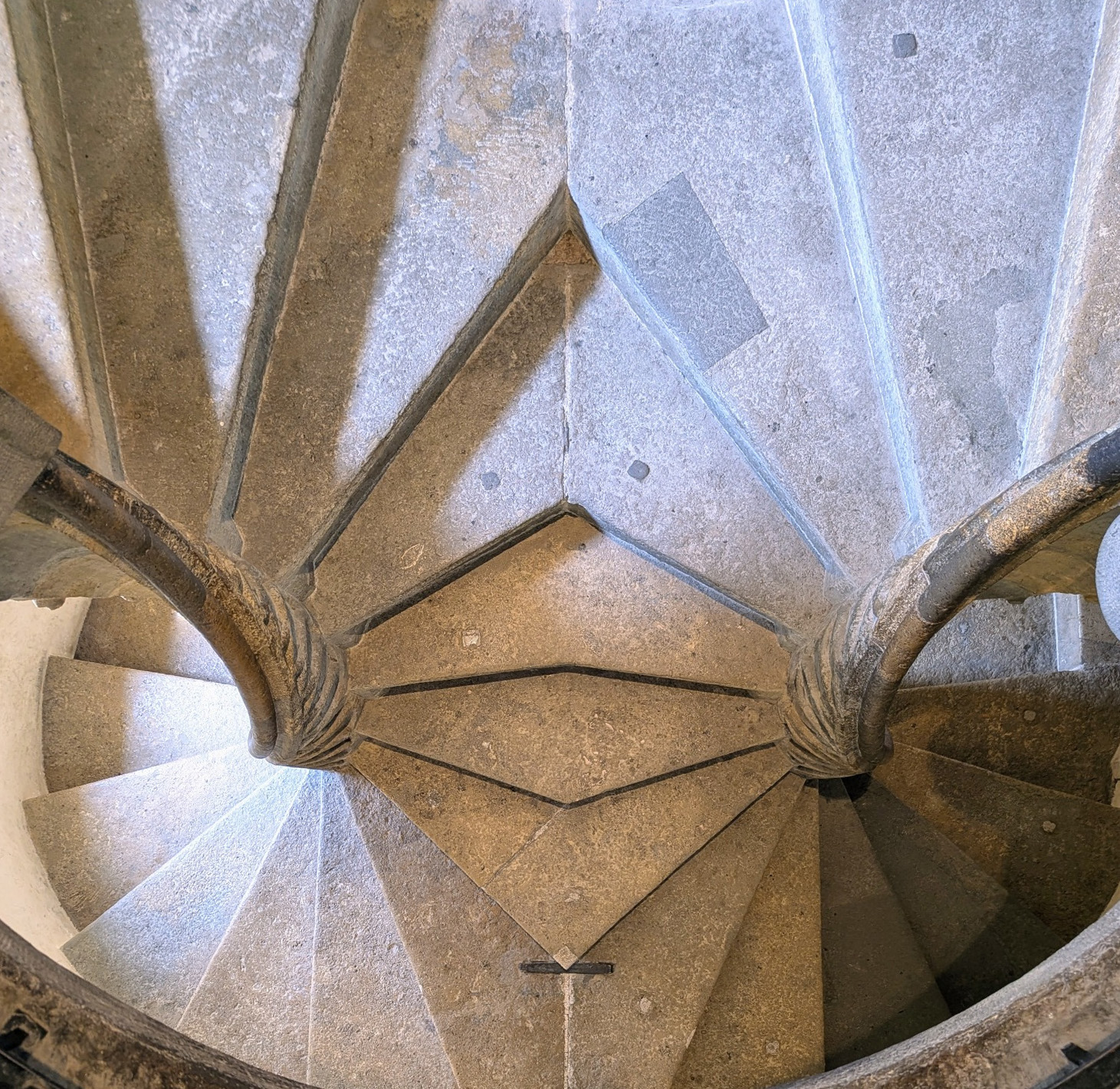}
  \caption{Looking down into a reunion zone. The two radial tread systems are
  organised around distinct spindles and become common through a polygonal
  sequence of shared stones. Photograph \copyright\ 2026 Santiago Schnell}
  \label{fig:overhead}
\end{figure}

The photographs make the geometry legible. Looking down from above
(Figure~\ref{fig:overhead}), one sees two radial fans organised around two
spindles. Near the middle, neither fan can continue as though the other were
absent. The treads become polygonal, align across the common zone, and then
recover their separate curvatures. The abstract intersections $P_+$ and $P_-$ are
not points one could stand upon; masonry gives them area, thickness, and
load-bearing form.

This is enough for the present essay. A full structural account would require a
surveyed geometry, material characterization, construction history, and an
explicit load-path analysis.

\section{Another double spiral staircase}
\begin{figure}[t]
  \centering
  \includegraphics[width=0.65\linewidth]{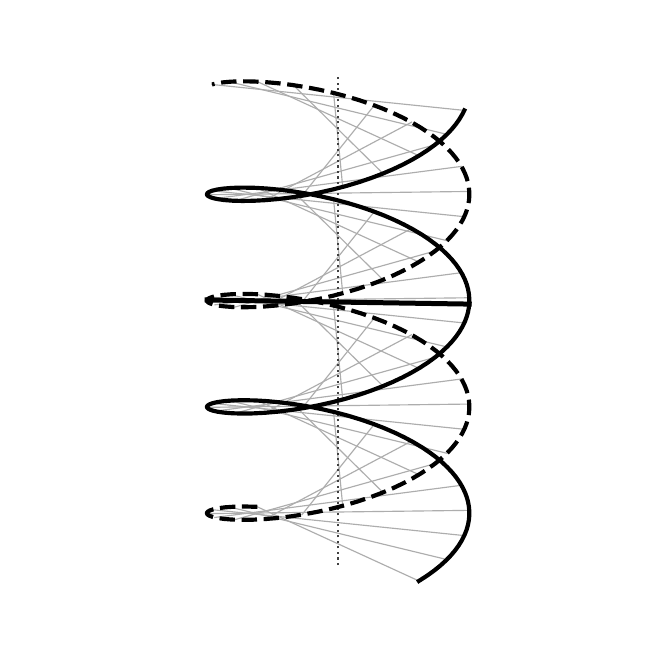}
  \caption{A truncated helicoid. The thin straight rulings rotate as their
  height changes; the two boundary curves are shown solid and dashed. They are
  phase-shifted, same-winding curves about one axis. Away from a reunion zone, a
  smoothed single flight of the Graz stair is locally modelled by an annular
  strip of such a surface; the twin system is not one helicoid}
  \label{fig:helicoid}
\end{figure}

There is, however, a different mathematical object that geometers also call a
``double spiral staircase'': the helicoid. For $c\ne0$, write
\begin{equation}
 H(u,v)=\bigl(u\cos v,\;u\sin v,\;cv\bigr),
 \qquad u,v\in\mathbb R.
 \label{eq:helicoid}
\end{equation}
For fixed $v$, the map $u\mapsto H(u,v)$ is a straight line. The helicoid is
therefore a ruled surface: a horizontal line rotates at a constant rate while
its height changes at a constant rate (Figure~\ref{fig:helicoid}).

The helicoid is a minimal surface. Indeed,
\[
 H_u=(\cos v,\sin v,0),\qquad
 H_v=(-u\sin v,u\cos v,c),
\]
so the first fundamental-form coefficients are
$E=1$, $F=0$, and $G=u^2+c^2$. A unit normal is proportional to
$(c\sin v,-c\cos v,u)$. The second fundamental-form coefficients are
\[
 e=g=0,\qquad f=\frac{-c}{\sqrt{u^{2}+c^{2}}},
\]
so the mean curvature is
\[
 \mathcal H=\frac{eG-2fF+gE}{2(EG-F^2)}=0.
\]
Thus the helicoid is stationary for the area functional. Its Gaussian 
curvature, by contrast, is
\[
 K=\frac{eg-f^{2}}{EG-F^{2}}=-\frac{c^{2}}{(u^{2}+c^{2})^{2}}<0 ,
\]
so the surface is saddle-shaped at every point: stationary for area, but nowhere
flat. Only in the limit $c\to0$ does it degenerate to the plane. Vanishing mean
curvature does not say that every compact piece is an area minimizer;
stability and minimizing properties require additional hypotheses
\citep{doCarmo1976}; for an illustrated survey of minimal surface forms, see
Almgren~\citep{Almgren1982}. Catalan~\citep{Catalan1842} proved in 1842 that the plane 
and the helicoid are the only ruled minimal surfaces.

Why is the helicoid called a double spiral staircase? Remove its axis and split
it into the regions $u>0$ and $u<0$. On $u=\rho>0$, the polar angle is
$\theta=v$ and $z=c\theta$. On $u=-\rho<0$, the polar angle is
$\theta=v+\pi$, so
\[
 z=c(\theta-\pi).
\]
Both sheets satisfy $dz/d\theta=c$. They spiral in the same
angular--vertical sense and are separated by a phase of $\pi$; they do not form
the oppositely handed, two-centre pair of the Graz stair. Colding and 
Minicozzi~\citep{ColdingMinicozzi2003,ColdingMinicozzi2006}
used precisely this double-staircase picture in their analysis of embedded
minimal disks: informally, regions where curvature concentrates are modelled,
after suitable rescaling, by two multi-valued sheets spiralling together around
a common axis.

Away from the reunion zones, however, a single smoothed flight has a natural
helicoidal idealization. Replacing its discrete treads by continuously rotating
radial segments gives an annular strip of a translated helicoid,
$r_{\mathrm{in}}\leq u\leq r_{\mathrm{out}}$. The two flights require distinct
axes and opposite signs of the pitch. This correspondence is local rather than
global: in each reunion zone, the masonry departs from two independent
helicoidal strips and resolves their overlap into common polygonal treads. The
Graz stair is therefore neither one helicoid nor simply the union of two
complete helicoids.

The comparison is summarized in Table~\ref{tab:comparison}.

\begin{table}[t]
\caption{Three constructions that are easily conflated. Their visual kinship is
real, but their axes, handedness, and modes of connection differ}
\label{tab:comparison}
\centering
\small
\begin{tabularx}{\linewidth}{@{}lYYY@{}}
\toprule
Object & Axes & Angular relation & How the two parts meet \\
\midrule
Graz twin stair (idealized) & Two & Opposite senses, with matched rise and phase & Centreline crossings become finite shared-tread zones \\
Helicoid split into two sheets & One & Same winding sense, phase shift $\pi$ & The sheets are disjoint after the axis is removed; their closures meet along the axis \\
Canonical B-DNA & One & Two same-handed helical strands, chemically antiparallel & The strands remain offset and connected by base pairing, not geometric crossing \\
\bottomrule
\end{tabularx}
\end{table}

The comparison with DNA is especially tempting to a mathematical biologist.
Canonical B-DNA consists of two strands in a right-handed double helix about one
axis \citep{WatsonCrick1953}. Its strands are chemically antiparallel, but that
is not the same as opposite geometric handedness. The distinction between the
Graz stair, the helicoid, and DNA is a small example of a large scientific rule:
analogy begins inquiry; it cannot replace identification.

\section{Reconciliation by design}
Correcting the geometry changes the moral of the staircase, and improves it.
Within the
symmetric recurrent model above, reunion is secured by $0<a<R$, equal radii,
equal rise rates, a synchronized phase, and a shared physical region in which
an intersection can become a tread. Alter the separation, radii, rise rates, or
phase, and the reunion may disappear. Reconciliation has
to be made possible.

I climbed the stair with more than a tourist's interest in that proposition. I
was trained first as a biologist and then took a doctorate in mathematics. I
have spent my working life since in the seam between them: biochemical
phenomena and dynamical systems, measurements and models, experiments and the
equations that would account for them. It has never felt like indecision. It has
felt like standing on the landing.

Mathematical biology is not biology decorated with equations, nor mathematics
using organisms as convenient examples. Its two practices meet only when a
common structure has been built: definitions refer to the same entities; units
and scales are compatible; experiments render parameters identifiable; models
preserve the biological meaning of their variables; and data are allowed to
correct the equations just as equations discipline measurement. The living 
world and the language of numbers belong, in the end, in the same council.

That, I think, is why the staircase would not let me go. Its flights remain
distinct. Each keeps its own centre and its own direction. Yet the mason gave
them overlapping ground, a matched ascent, and stones cut to bear both paths.
The result is not an infinity trap but a repeated invitation: separate, meet,
and separate again without pretending that difference has vanished. It is a
good thing to be reminded of halfway up a tower, between sessions.

\section{For the mathematical tourist}
The entrance to the \emph{Doppelwendeltreppe} is at Hofgasse~15, inside the courtyard of
the Grazer Burg. From Graz Cathedral, walk a short distance north along Hofgasse,
enter the Burg courtyard, and look for the narrow doorway near the arch bearing
the A.E.I.O.U. inscription. At the time of writing, admission is free; current
opening hours should be checked on the official visitor page
\citep{GrazTourismus}.


\backmatter

\section*{Acknowledgments}
I thank my wife, Mariana Schnell, for discovering the staircase, leading me to it,
and sharing the photographic record of our visit. I also thank the organisers and
participants of ECMTB 2026, and the leadership of the European Society for
Mathematical and Theoretical Biology and the Society for Mathematical Biology,
whose meeting supplied the occasion for this tourist visit.

\section*{Declarations}

\noindent\textbf{Funding} This work was supported by Dartmouth College.

\noindent\textbf{Competing interests} The author has no relevant financial or non-financial interests to disclose.

\noindent\textbf{Ethics approval and consent to participate} Not applicable.

\noindent\textbf{Consent for publication} Not applicable.

\noindent\textbf{Data availability} Data sharing is not applicable to this article because no datasets were generated or analysed.

\noindent\textbf{Materials availability} Not applicable.

\noindent\textbf{Code availability}
No computational results are reported, and no research code is required to
reproduce the mathematical arguments or conclusions of this article.
Figures~\ref{fig:plan} and~\ref{fig:helicoid} are illustrative renderings
of equations stated explicitly in the text.

\noindent\textbf{Author contributions} Santiago Schnell conceived the article, developed and checked the mathematical model, selected the illustrations, and wrote and revised the manuscript with the support of an AI editing tool.

\bibliography{references}

\end{document}